\documentclass[12pt]{article}
\usepackage{amsmath}
\usepackage{amssymb}
\usepackage[]{graphicx}
\usepackage{epsfig, float}
\usepackage[cp1251]{inputenc}
\usepackage[english]{babel}
\usepackage{color}

\newtheorem{proposition}{Proposition}

\newcommand{\be}{\begin{equation}}
\newcommand{\ee}{\end{equation}}
\newcommand{\bea}{\begin{eqnarray}}
\newcommand{\eea}{\end{eqnarray}}

\newcommand{\ep}{\varepsilon}

\newtheorem{theorem}{Theorem}
\newtheorem{remark}{Remark}

\newtheorem{lemma}{Lemma}

 \newcounter{myc}

\begin{document}
 \date{}

\title{On strong duality in linear copositive programming
\thanks{This   work was partially supported state research program "Convergence"(Republic Belarus), Task 1.3.01 and by  Portuguese funds through CIDMA - Center for Research and Development in Mathematics and Applications, and FCT -
 Portuguese Foundation for Science and Technology, within the project  UID/MAT/04106/2019.}}

\author{Kostyukova O.I.\thanks{Institute of Mathematics, National Academy of Sciences of Belarus, Surganov str. 11, 220072, Minsk,
 Belarus  ({\tt kostyukova@im.bas-net.by}).} \and Tchemisova T.V.\thanks{Mathematical Department, University of Aveiro, Campus
Universitario Santiago, 3810-193, Aveiro, Portugal ({\tt tatiana@mat.ua.pt}).}}
\maketitle

\begin{abstract}

The paper is dedicated to the study of strong duality  for a problem of linear copositive programming.
 Based on the recently introduced concept of the set of normalized immobile indices,
an extended dual problem is  deduced. The dual problem satisfies the strong duality relations
and does not {require} any additional regularity assumptions such as constraint qualifications.
The main difference with {the} previously obtained results consists in the fact that now the extended dual problem
  {uses neither} the immobile indices themselves {nor the explicit information about the convex hull of these indices}.

 The strong duality formulations presented in the paper have similar structure and properties
as that proposed in the works of  M. Ramana, L. Tuncel, and H. Wolkovicz, for semidefinite programming,
but are obtained using different techniques.

\end{abstract}

\\


\section{Introduction}\label{Intro}
Linear {\it Copositive Programming} problems can be considered as linear programs over the convex cone of so-called {\it copositive matrices} (i.e. matrices which are positive semi-defined on the non-negative orthant). Copositive problems form a special class of conic optimization problems and have many important applications, including $\mathcal{NP}$ -hard problems. For the references on applications  of Copositive Programming see e.g. \cite{Bomze 2012, Dur 2010} and others.

 Linear copositive problems are closely related to that of linear {\it Semi-Infinite Programming} (SIP) and
{\it Semidefinite Programming} (SDP). Linear copositive and semidefinite problems are particular cases of SIP problems but Linear Copositive Programming deals with  more challenging and less studied  problems   than that of SDP. The literature on theory and methods of SIP, SDP and Linear Copositive Programming is very interesting and quite large, we refer the interested readers to \cite{Bomze 2012,Bon,Dur 2010,Kort, Li,Luo, Ramana,Ramana-W, Weber2}, and the references therein.

{Although} the concepts of  co-positivity and complete positivity were originally formulated in 1952
in the paper by T. Motzkin \cite{Motzkin}, an active research in theory and methods of Linear Copositive Programming has
 begun only in the recent decades in papers of I. Bomze, M, D\"{u}r,  E. de Klerk, and others (see \cite{Bomze 2000,  Bomze 2012, Dick}).

Optimality conditions and  the associated  duality relationships are among the central topics
of convex optimization and the   importance of their  study is well recognized (see e.g. \cite{Bon, Shapiro2}, and the references therein).
 Optimality conditions are a crucial issue in the study of any optimization problem since they
allow not only to test the optimality of a given feasible  solution,  but also to develop efficient numerical methods.
 As it is mentioned in \cite{Luo}, the duality plays a central role in detecting infeasibility,
lower-bounding of the optimal objective value, as well as in design and analysis of iterative algorithms.

 Often the studies on optimality conditions and duality for finite and (semi-) infinite programming  use certain regularity assumptions,
 so-called {\it constraint qualifications} (CQ). Such assumptions permit to guarantee in some particular cases a {\it strong}
(or {\it zero-gap}) duality which means that the optimal values of the primal and dual objective functions are equal and, hence,
the difference between these values (the duality gap) vanishes.

It is a known fact that in  Linear Programming (LP), the   strong duality is guaranteed
without any CQ  (\cite{Bon}). The duality results for LP can be generalized to some particular classes of optimization problems.
Several attempts were done to obtain CQ-free optimality and  strong   duality results for different classes of convex SIP problems (see e.g. \cite{Jeyak,
KT-JOTA, KT-Opt}).

In \cite{Ramana-W, Tuncel}, a CQ-free duality theory for conic optimization  was developed
in terms of so-called {\it minimal cone}.
Being quite general, this theory has one disadvantage in terms of its application, namely, it is very abstract.
In a number of publications, various explicit dual
 formulations were obtained {by applying} this theory to SDP   problems  (see e.g. \cite{Kort, Ramana, Shapiro2, Tuncel}) and {other}
optimization problems over symmetric (i.e., self-dual and homogeneous)
 cones (see \cite{Terlaki}).  As it was mentioned in \cite{PP,Terlaki}, finding a broader family
of conic problems for which such  explicit dual  formulations  are possible, is an open problem.

Linear copositive problems   belong to a wider and  more complex class of linear conic problems than that of SDP, namely, to the
class  of  optimization problems over  cones  of  copositive and$\setminus$or  completely positive matrices that are
neither self-dual nor homogeneous (see \cite{GOWDA}). The duality theory for these problems is not   well studied yet.
It is worth to mention that almost all duality results and optimality conditions  for Linear Copositive Programming are
 formulated under the Slater CQ (\cite{AH 2013,Bomze 2000}).

In our papers \cite{KT-JOTA, KT-Opt, KT-JMS}, and others, we developed  a new approach to optimality in  SIP and SDP.
This approach is based on the notion of {\it immobile indices} of constraints of an optimization problem, which refers to the
 indices of the constraints that are active for all feasible solutions.

In \cite{KT-SetValued, KT-Opt2018}, we have applied our approach to
problems of Linear Copositive Programming  and
 successfully obtained  new explicit  CQ-free optimality conditions and  strong duality results. It is essential that to  formulate  our results, we
 used either the immobile indices (\cite{KT-Opt2018}), or the vertices of the convex hull of the {\it normalized immobile index set} (\cite{KT-SetValued}).

In this paper, we further develop our approach to linear copositive problems and use it to obtain a new  dual  problem which
 we refer  to  here as  the {\it extended dual problem}.  As well as the {\it regularized dual problem}  from \cite{KT-SetValued},
 the extended dual one is  constructed using  the notion and properties of the normalized immobile indices, but in the formulation of this problem
 neither  these indices nor the vertices of the convex hull of the corresponding index set are present.
 This   permitted us to formulate the extended dual  problem  for linear copositive problem  in  an explicit   form and  avoid  the  use of additional procedures for finding the immobile indices. The new extended dual problem satisfies
the strong duality relations without any CQ.

  One interesting property of the obtained results consists in the fact that the new dual formulations
	for Linear Copositive Programming are closely related to the dual problems proposed in  \cite{Ramana, Ramana-W, Tuncel} for SDP.
	This relation not only confirms the already known deep connection between copositive and semidefinite problems     but,
	taking into account the impact of the duality results of M. Ramana  et al. in relation to SDP,
	permits one to expect that the proposed here  duality  results are also very promising in  Linear Copositive Programming.

 It is worth to mention here that at present, with exception of \cite{KT-SetValued}, there are no explicit strong duality formulations without CQs for
 Linear Copositive Progamming. 	All the results presented in the paper  are original and cannot be obtained as a direct extension of any previous results.

\vspace{2mm}

The paper is organized as follows. Section 1  hosts Introduction. In section 2,
given a linear copositive problem, we formulate the corresponding normalized immobile index set and establish
 some new properties of this set. An extended dual problem is formulated in section 3. We  prove here  that  the {strong} duality property is satisfied.
In section 4, we compare
the obtained duality results   with that  presented in \cite{Ramana, Ramana-W, Tuncel} for SDP.  It is shown  that the
 compared dual formulations  for Linear Copositive Programming and SDP are similar, being both  CQ-free and providing strong duality.
 Although  these dual formulations were obtained using different techniques, they  almost coincide being applied to the class of
 linear SDP problems. The final section 5 contains some conclusions.

\section{Linear copositive  programming  problem }\label{optimality}

Here and in what follows, we use  the following notations.
Given an integer $p>1$, $\mathbb R^{ p}_+$ denotes the set of all $ p$ vectors with non-negative components,
 ${\cal S}(p)$ stays for the space of real symmetric $p\times p$ matrices,
${\cal P}(p)$  for the cone of symmetric positive semidefinite $p\times p$ matrices, and
$\mathcal{COP}^p$ for the cone of symmetric copositive $p\times p$ matrices
\begin{equation}\label{CopositiveT}\mathcal{COP}^p:=\{D\in {\cal S}(p):t^{\top}Dt\geq 0 \ \forall t \in \mathbb
R^p_+\}.\end{equation}

The space $\mathcal S(p)$ is considered here as a vector space with the trace inner product: $$A\bullet B:={\rm trace}\, (AB), \ \mbox{ for } A,B \in \mathcal S(p).$$

Consider a linear Copositive Programming  problem in the form
\begin{equation} \displaystyle\min_{x} c^{\top}x \;\;\;
\mbox{s.t. } {\cal A}(x)
\in \mathcal{COP}^p,\label{SDP-2}\end{equation}
 where the decision variable is $n-$vector $x=(x_1,...,x_n)^{\top}$  and the constraints matrix function $\mathcal A(x)$ is defined as
\begin{equation*}\label{A} \mathcal A(x):=\displaystyle\sum_{i=1}^{n}A_ix_i+A_0,\end{equation*}
matrices $A_i\in \mathcal S(p), i=0,1,\dots,n,$  and  vector $c \in \mathbb R^n $   are given.
Problem (\ref{SDP-2}) can be rewritten as follows:
\begin{equation} \displaystyle\min_{x} c^{\top}x \;\;\;
\mbox{s.t. } t^{\top}{\cal A}(x)t
\geq  0 \;\; \forall t \in \mathbb R^p_+.\label{SDP-1}\end{equation}

It is well known that the  copositive problem (\ref{SDP-1}) is equivalent to the following convex  SIP problem:
\begin{equation} \displaystyle\min_{x} c^{\top}x \;\;\;
\mbox{s.t. } t^{\top}{\cal A}(x)t
\geq  0 \;\; \forall t \in T,\label{SIP}\end{equation}
 with a $p$ - dimensional compact index set in the form of a simplex
 \begin{equation}\label{SetT}T=\{t \in \mathbb R^p_+:\mathbf{e}^{\top}t=1\},\end{equation}
where $\mathbf{e}=(1,1,...,1)^{\top}\in \mathbb R^p,$  $t=(t_k,k\in P)^{\top}, $ $ P=\{1,2,...,p\}.$

Denote by $X$ the set of feasible solutions of problems (\ref{SDP-2}) -  (\ref{SIP}),
$$X:=\{x\in \mathbb R^n:t^{\top}{\cal A}(x)t\geq 0 \;\; \forall  t \in \mathbb R^p_+\}
=\{x\in \mathbb R^n:t^{\top}{\cal A}(x)t\geq 0 \;\; \forall  t \in T\}.$$

Evidently,  the set $X$ is convex.

 According to the definition (see e.g. \cite{KT-Opt2018}), the constraints of the SIP problem (\ref{SIP}) satisfy the Slater condition if
\begin{equation}\label{Slait}
\exists \; \bar{x} \in \mathbb R^n \; \mbox{ such that } \; t^{\top}{\cal A}(\bar x)t >0 \  \forall t\in T,
\end{equation}
 and the constraints of the copositive problem (\ref{SDP-2}) satisfy the Slater condition if
\begin{equation}\label{Slait1}
\exists \; \bar{x} \in \mathbb R^n \; \mbox{ such that } \; {\cal A}(\bar x) \in {\rm int }\, \mathcal{COP}^p=
\{D\in \mathcal{S}(p): t^{\top}Dt>0 \; \forall t \in \mathbb R^p_+,\; t \not=0\}.
\end{equation}
Here ${\rm int  }\,D$ stays for the interior of a set $D$.

Evidently, problems (\ref{SDP-2}), (\ref{SDP-1}), and (\ref{SIP})  satisfy or not the Slater condition simultaneously.

Following \cite{KT-JOTA, KT-SetValued}, define the sets of immobile indices $ T_{im}$  and $ R_{im}$ in  problems (\ref{SIP}) and  (\ref{SDP-1}), respectively:

\vspace{-3mm}

\begin{equation*}\label{T*} T_{im}:=\{t \in   T : t^{\top} \mathcal A (x)t=0 \ \  \forall x\in {X}\}\end{equation*}
and
\begin{equation*}\label{hatT*}  R_{im}:=\{t \in  \mathbb R^p_+ : t^{\top} \mathcal A (x)t=0 \ \ \forall x\in {X}\}.\end{equation*}

It is evident that the aforementioned sets  are interrelated:
$$  R_{im}=\{t \in \mathbb R^p: t=\alpha \tau , \;\; \tau \in T_{im}, \; \alpha\geq 0\} \ \mbox{and}\
T_{im}=\{t \in R_{im}: \mathbf{e}^{\top}t=1\}.$$

From the latter  relations, we conclude that the set $T_{im}$, the {\it immobile index set} for problem (\ref{SIP}), can be considered as a {\it
normalized immobile index set} for problem (\ref{SDP-1}).  In what follows, we will use mainly the set $T_{im}$, taking into account its relationship with the set $R_{im}$.

 The following proposition is an evident corollary of Proposition 1 from \cite{KT-SetValued}.

 \vspace{-1mm}

\begin{proposition}\label{P3} Given a linear copositive problem in the form (\ref{SDP-1}), the Slater condition (\ref{Slait}) is equivalent to the emptiness of the normalized
immobile index set $T_{im}.$
\end{proposition}

It is evident that $T_{im}=\emptyset$ if and only if $R_{im}=\{\mathbf{0}\}$.

\begin{proposition}\label{Pnew1} Given a  linear copositive problem (\ref{SDP-1}), let $\{\tau(i), i \in I\}$ be
some set consisting of  immobile indices of this problem. Then for any $x\in X,$ the following inequalities take a place:
\be {\cal A}(x)\tau(i) \geq 0, \;i\in I. \label{3**}\ee

\end{proposition}

The proof of the proposition follows from the definition of immobile indices and Lemma 2.6 from \cite{Bam}.

\vspace{5mm}

\begin{proposition}\label{Pnew2} Given an index set $\{\tau(i)\in T, i \in I\},$  the inequalities
\be {\cal A}(x)\tau(i)\geq 0, \;i\in I, \label{03**}\ee
imply the inequalities
\be t^{\top}{\cal A}(x)t\geq 0\; \forall t \in {\rm conv} \{\tau(i), i \in I\}.\label{19-1}\ee
\end{proposition}
Here ${\rm conv} S$ denotes the convex hull of a given set $S$.

The proof of the proposition is evident.

\vspace{5mm}

Let $\{\tau(i), i \in I\}\subset T_{im}$ be a nonempty subset of the set of normalized immobile indices in  problem (\ref{SDP-1}).
For this subset and for any $\varepsilon>0$, denote

\vspace{-4mm}

\begin{equation}\label{Tep} T(\varepsilon):=T(\varepsilon,\tau(i), i \in I):=\{t \in T:\ \rho(t,{\rm conv}\ \{\tau(i), i \in
I\})\geq {\ep}\},\end{equation}
 \begin{equation}\label{Tep1}\widehat T(\ep):=\widehat T(\varepsilon,\tau(i), i \in I):=\{t \in T: \ \rho(t,{\rm conv}\
 \{\tau(i), i \in  I\})\leq \ep\},\end{equation}

where $\rho(l,B)=\min\limits_{\tau \in B}||l-\tau||$ is the distance between a vector $l$ and a set $B$  associated with the norm
$||a||=\sqrt{a^{\top}a}$ in the vector space $\mathbb R^p$.
Consider the sets
\be\mathcal X=\{x \in \mathbb R^n: {\cal A}(x)\tau(i)\geq 0, \;i\in I\},\ \
 \mathcal X(\ep)=\{z \in \mathcal X: t^{\top}{\cal A}(z)t\geq 0, \;\forall  t \in
 T(\varepsilon)\}.\label{3}\ee

The  following lemma  is a generalization of Lemma 2 from \cite{KT-SetValued}.
\begin{lemma}\label{L0} Let $\{\tau(i), i \in I\}$ be a subset of the set of normalized immobile indices in  problem (\ref{SDP-1}).  Then there exists
$\ep_0>0$ such that
$\mathcal X(\ep_0)=X$,  the set ${\mathcal X(\ep)}$ being defined in (\ref{3}) with the set $T(\ep)$  as in (\ref{Tep}).

\end{lemma}

\textbf{Proof. }
 It follows from Proposition \ref{Pnew1}  that $X\subset \mathcal X(\ep)$ for all $ \ep>0.$
To finalize the proof, it is enough to show that there exists $\ep_0>0$ such that $\mathcal X(\ep_0)\subset X.$
Suppose the contrary. Then  for each $\ep >0$  there exist $ z(\ep) \in \mathcal X(\ep)$ such that
\be (t(\ep))^{\top}{\cal A}(z(\ep))t(\ep)<0,\label{33-4}\ee
where $t(\ep)$ is an optimal solution of the problem
\begin{equation}\label{LLP-ep}
\min\limits_{t\in T}t^{\top}{\cal A}(z(\ep))t.
\end{equation}
 Since, by construction (see Proposition \ref {Pnew2} and (\ref{3})), it holds $$t^{\top}  {\cal A}(z(\ep))t\geq 0 \; \forall t \in
 T(\varepsilon)\cup {\rm conv} \ \{\tau(i),i \in I\},$$

then
 $t(\ep)\in \widehat T(\varepsilon)\setminus{\rm conv}\ \{\tau(i),i \in I\}$ for all $\ep >0$.
Hence there exists $t^*:=\lim\limits_{\ep\to 0}t(\ep),$  $ t^*\in {\rm conv}\ \{\tau(i),i \in I\}.$

For  $\ep>0$, let us consider the vector $ l(\ep):=t(\ep)-t^*.$
By construction, $e^{\top} l(\ep)=0$.

It is evident that  there exists a  sufficiently small  $\bar \ep>0$ such that  for $k\in P$, the following conditions hold:
\begin{center}
if $t^*_k=0,$ then $ l_k(\bar \ep)=t_k(\bar \ep)\geq 0 \ $  and if $ \ t_k(\bar \ep)= 0,$  then $t^*_k=0$ and $ l_k(\bar \ep)=0.$
\end{center}

Consequently, the direction $l:= l(\bar \ep)$ is feasible for $t^*$ and $t(\bar \ep)$ in the set $T$.  Hence
there exists $\gamma_0>1$ such that
$t^*+\gamma  l=t^*+\gamma(t(\bar \ep)-t^*)\geq 0, \; {\mathbf e}^{\top}{(t^*+\gamma  l)}=1\;\; \forall \gamma\in [0,\gamma_0].$

Define the function
\begin{equation*}\begin{split}w(\gamma)&:={(t^*+\gamma  l)^\top}{\cal A}(z(\bar \ep)){(t^*+\gamma  l)}\\&=(t^*)^{\top}{\cal A}(z(\bar \ep))t^*+2\gamma
l^{\;\top}{\cal A}(z(\bar \ep))t^*+
\gamma^2 l^{\;\top}{\cal A}(z(\bar \ep))\bar l=a\gamma^2+2b\gamma+c,\end{split}\end{equation*}

\vspace{-2mm}

where $c:= (t^*)^{\top}{\cal A}(z(\bar \ep))t^*,$ $b:= l^{\; \top}{\cal A}(z(\bar \ep))t^*,$
$a:= l^{\; \top}{\cal A}(z(\bar \ep)) l.$

\vspace{2mm}

By construction, for $\gamma^*:=1$ we have
$w(\gamma^*)=t^{\top}(\bar \ep){\cal A}(z(\bar \ep))t(\bar \ep)$ and it is the optimal value of the cost function of the problem (\ref{LLP-ep}) with $\ep =\bar
\ep.$
Hence
$$ w(\gamma^*)=\min\limits_{\gamma\in [0,\gamma_0]}w(\gamma)= \min\limits_{\gamma\in [0,\gamma_0]}(a\gamma^2+2b\gamma+c).$$
Since $\gamma^*=1\in (0,\gamma_0)$ in the formula above, then $2a\gamma^*+2b=2a+2b=0$.
Therefore
 $-b=a$, which can be rewritten in the form
$- l^{\; \top}{\cal A}(z(\bar \ep))t^*= l^{\; \top}{\cal A}(z(\bar \ep)) l,$
{wherefrom} we get
\be (t(\bar \ep))^{\top}{\cal A}(z(\bar \ep))t^*=(t(\bar \ep))^{\top}{\cal A}(z(\bar \ep))t(\bar \ep).\label{33-7}\ee
Since $t^*\in {\rm conv}\ \{\tau(i),i \in I\}$, then
$t^*=\sum\limits_{i\in I}\beta_i \tau(i),\; \;\sum\limits_{i\in I}\beta_i=1,\;  \beta_i\geq 0, \; i\in I.$

Consequently, taking into account that $z(\bar \ep) \in  {\mathcal X(\bar \varepsilon)} $  and $t(\bar \ep)\geq 0$, we
have
$$  (t(\bar \ep))^{\top}{\cal A}(z(\bar \ep))t^*=\sum\limits_{i\in I}\beta_i (t(\bar \ep))^{\top}{\cal A}(z(\bar \ep))\tau(i)\geq 0.$$
But this inequality and inequality (\ref{33-4}) contradict equality (\ref{33-7}).
The lemma is proved. $\Box$
\vspace{5mm}

 It should be noticed that the lemma above can be considered as a generalization of Lemma 2 from \cite{KT-SetValued} since
it is formulated for   an {\bf arbitrary } subset $\{\tau(i), i \in I\}$ of $T_{im}$,
while in  Lemma 2 from \cite{KT-SetValued} we considered the {\bf fixed} subset of $T_{im}$, namely, the set of vertices
of ${\rm conv} T_{im}$.

\section{An extended dual problem for Linear Copositive\\ Programming}\label{ExDual}

In this section, we will  formulate an extended dual problem for  problem {(\ref{SDP-2})}.

Given an arbitrary  cone ${\cal K }\in {\cal S}(p)$, the corresponding dual cone ${\cal K}^*$ is defined as
$${\cal K }^*:=\{A\in {\cal S}(p): \ A\bullet D\geq 0\;\; \forall  D\in {\cal K}\}.$$

It is known that  the cone of symmetric positive  semidefinite matrices ${\cal P}(p)$ is self-dual, i.e. ${\cal P}^*(p)={\cal P}(p)$ but the cone of
symmetric  copositive matrices ${\mathcal{COP}^p}$ defined in (\ref{CopositiveT}), is
not.

 It can be shown (see e.g. \cite{Abraham}) that  for {the} cone ${\mathcal{COP}^p}$,   its
dual cone ${\mathcal{CP}^p}$ is the cone of so-called completely positive matrices,
${\mathcal{CP}^p:}= {\rm conv}\{x x^{\top}: x \in {\mathbb R}^p_+\},$  and  ${\mathcal{CP}^p}\subset
{\mathcal{COP}^p}.$

For a given finite integer $m_0\geq 0$, consider the following problem:
\begin{equation}\label{dual}\begin{split}
&  \max \;\; -(U+W_{m_0})\bullet A_0,\\
\mbox{s.t. \qquad }&(U_m+W_{m-1})\bullet A_j=0,\; j=0,1,...,n,\; m=1,...,m_0;\\
& (U+W_{m_0})\bullet A_j=c_j, \; j=1,2,...,n;\\
& U\in \mathcal{CP}^p,\;\; W_0=\mathbb O_{p},\\
& \left(\begin{array}{cc}
U_m& W_m\cr
W_m^{\top}& D_m
\end{array}\right)\in \mathcal{CP}^{2p},\; m=1,...,m_0,\end{split}\end{equation}
where $U_m\in {\cal S}(p),\; D_m\in {\cal S}(p),$  $W_m\in \mathbb R^{p\times p},  m=1,...,m_0,$ and $\mathbb O_{p}$ stays for the
 $p\times p$ null matrix.

  Notice that in the case $m_0=0$, we consider that the index set $\{1,...,m_0\}$ is empty and
the  constraints $(U_m+W_{m-1})\bullet A_j=0,\; j=0,1,...,n,$ $\left(\begin{array}{cc}
U_m& W_m\cr
W_m^{\top}& D_m
\end{array}\right)\in \mathcal{CP}^{2p},$ $\; m=1,...,m_0,$  are absent in problem
(\ref{dual}).  Hence,   for $m_0=0$,   problem (\ref{dual}) takes the form
\begin{equation}\label{dual0}\begin{split}
& \max  \; -U\bullet A_0,\\
\mbox{s.t.} \qquad & U\bullet A_j=c_j, \; j=1,2,...,n;\;\;
 U\in \mathcal{CP}^p.\end{split}\end{equation}

\begin{lemma}\label{17-L}[Weak duality]
Let $x\in X$ be a feasible solution of the primal linear copositive problem (\ref{SDP-2}) and
\be (U_m,\,W_m,\, D_m, \; m=1,...,m_0;\;\; U)\label{dGam}\ee
be a feasible solution of problem (\ref{dual}). Then the following inequality holds:
\be c^{\top}x\geq -(U+W_{m_0})\bullet A_0.\label{17-3}\ee
\end{lemma}
{\bf Proof.}  For $m=1,...,m_0, $ it follows from the condition $ \left(\begin{array}{cc}
U_m& W_m\cr
W_m^{\top}& D_m
\end{array}\right)\in \mathcal{CP}^{2p}$ that there exists a matrix $B_m$
with non-negative elements in the form
$$B_m=\left(\begin{array}{c}
V_m\cr
L_m
\end{array}\right)\in\mathbb R^{2p\times k(m)},$$
 such that
$$\left(\begin{array}{cc}
U_m& W_m\cr
W_m^{\top}& D_m
\end{array}\right)=B_mB^{\top}_m=\left(\begin{array}{c}
V_m\cr
L_m
\end{array}\right)(V^{\top}_m\; \; L^{\top}_m).$$
{The matrix $B_m$ above is composed by the blocks containing some matrices
\bea & V_m=(\tau^m(i), i \in I_m), \; L_m=(\lambda^m(i), i \in I_m), \label{Vm}\\
\mbox{where}&\tau^m(i)\in \mathbb R^p_+, \; \lambda^m(i)\in \mathbb R^p_+, i \in I_m, \; k(m):=|I_m|,\nonumber\eea}
 Hence, for  $m=1,...,m_0,$ the matrices $U_m,\; W_m,\; D_m$ in (\ref{dual}) admit  representations
\begin{equation}U_m=V_mV^{\top}_m,\;\; W_m=V_mL^{\top}_m,\; D_m=L_m
L^{\top}_m.\label{REP}\end{equation}

Consider the first group of constraints of the dual problem (\ref{dual}):
$U_1\bullet A_j=0,\; j=0,1,...,n.$
Due to  (\ref{Vm}) and (\ref{REP}), these constraints can be rewritten in the form
\be \sum\limits_{i\in I_1}(\tau^1(i))^{\top}A_j\tau^1(i)=0,\; j=0,1,...,n.\label{17-4}\ee
 It follows from (\ref{17-4}) that for any $x\in \mathbb R^n,$ we have
\be \sum\limits_{i\in I_1}(\tau^1(i))^{\top}{\cal A}(x)(\tau^1(i))=0.\label{17-5}\ee
Taking into account that for any $x\in X$, the inequalities
\be t^{\top}{\cal A}(x)t\geq 0 \; \forall t \in \mathbb R^p_+,\label{17-17}\ee
should  be fulfilled,   equality (\ref{17-5}) implies
$(\tau^1(i))^{\top}{\cal A}(x)\tau^1(i)=0,$ $ i \in I_1, \;\; \forall x \in X.$

Thus one can conclude that  $\tau^1(i)\in T_{im}, i \in I_1,$
and, consequently (see Proposition \ref{Pnew1}),
\be {\cal A}(x)\tau^1(i)\geq 0,\;\; i \in I_1,\  \forall x \in X.\label{17-5.1}\ee

Suppose that for some $m\geq 1$, it was shown that
\be {\cal A}(x)\tau^m(i)\geq 0,\;\; i \in I_m,\ \forall  x \in X.\label{17-8}\ee
Due to (\ref{Vm})  and (\ref{REP}), the constraints $(U_{m+1}+W_{m})\bullet A_j=0, $ $j=0,1,...,n,$ of problem (\ref{dual})
 can be rewritten as follows:
$$ \sum\limits_{i\in I_{m+1}}(\tau^{m+1}(i))^{\top}A_j\tau^{m+1}(i)+
\sum\limits_{i\in I_{m}}(\lambda^{m}(i))^{\top}A_j\tau^{m}(i)=0,\; j=0,1,...,n.$$
It follows from the {latter}  equalities  that for any $x\in \mathbb R^n,$ we have
\be \sum\limits_{i\in I_{m+1}}(\tau^{m+1}(i))^{\top}{\cal A}(x)\tau^{m+1}(i)+
\sum\limits_{i\in I_{m}}(\lambda^{m}(i))^{\top}{\cal A}(x)\tau^{m}(i)=0.\label{17-9}\ee
By the hypothesis above,  inequalities (\ref{17-8}) are satisfied. Then,
 taking into account that $\lambda^{m}(i)\in \mathbb R^p_+,$ $i \in I_m,$ and  for any $x\in X$,  inequalities (\ref{17-17}) hold,  we conclude from
  (\ref{17-9}) that
$$(\tau^{m+1}(i))^{\top}{\cal A}(x)\tau^{m+1}(i)=0,\;\; i \in I_{m+1},\;\; \;(\lambda^{m}(i))^{\top}{\cal A}(x)\tau^{m}(i)=0,\;  i\in I_m,
 \; \forall  x \in X.$$
Hence, $\tau^{m+1}(i)\in T_{im}, i \in I_{m+1},$  and, according to Proposition~\ref{Pnew1}, it holds
$$ {\cal A}(x)\tau^{m+1}(i)\geq 0,\;\; i \in I_{m+1},\; \forall  x \in X.$$

Now, replace $m$ by $m+1$ and repeat the considerations for all  $m<m_0.$

Let $m=m_0$. In this case,   relations (\ref{17-8}) have the form
\be {\cal A}(x)\tau^{m_0}(i)\geq 0,\;\; i \in I_{m_0},\;\; \forall  x \in X,\label{17-10.1}\ee
and for $U=\sum\limits_{i\in I}\tau(i)\tau^{\top}(i), \; \tau(i)\in \mathbb R^p_+,\; i \in I,$ the constraints $$(U+W_{m_0})\bullet A_j=c_j,\;\; j=1,...,n,$$
 of problem (\ref{dual})  can be represented as  follows:
\be \sum\limits_{i\in I}(\tau(i))^{\top}A_j\tau(i)+
\sum\limits_{i\in I_{m_0}}(\lambda^{m_0}(i))^{\top}A_j\tau^{m_0}(i)=c_j, \; j=1,...,n.\label{17-10}\ee

Then, evidently,
\begin{equation}\begin{split} \sum\limits_{j=1}^n c_jx_j= &\sum\limits_{i\in I}(\tau(i))^{\top}{\cal A}(x)\tau(i)+
\sum\limits_{i\in I_{m_0}}(\lambda^{m_0}(i))^{\top}{\cal A}(x)\tau^{m_0}(i)\\&-\Big(\sum\limits_{i\in I}(\tau(i))^{\top}A_0\tau(i)+
\sum\limits_{i\in I_{m_0}}(\lambda^{m_0}(i))^{\top}A_0\tau^{m_0}(i)\Big)\\ &=
\sum\limits_{i\in I}(\tau(i))^{\top}{\cal A}(x)\tau(i)+
\sum\limits_{i\in I_{m_0}}(\lambda^{m_0}(i))^{\top}{\cal A}(x)\tau^{m_0}(i)-(U+W_{m_0})\bullet A_0.\label{17-11-z}\end{split}\end{equation}
\normalsize{}
From (\ref{17-17}) and (\ref{17-10.1}) we conclude that
$$(\tau(i))^{\top}{\cal A}(x)\tau(i)\geq 0, \; i \in I;\;\; (\lambda^{m_0}(i))^{\top}{\cal A}(x)\tau^{m_0}(i)\geq 0,\; i\in I_{m_0},\;\; \forall  x \in X.$$
These  inequalities together with equality (\ref{17-11-z}) imply (\ref{17-3}). The lemma is proved.
$\qquad \Box$

\begin{lemma}\label{17-L-2}[Strong duality] Let problem (\ref{SDP-2}) have an optimal solution. Then there exist a number $0\leq m_0<\infty$ and a feasible
solution
\be \Bigl(U^0_m,\,W^0_m,\, D^0_m,\; m=1,...,m_0;\; U^0\Bigl)\label{G0}\ee

of problem (\ref{dual}) such that for any  optimal solution $x^0$ of problem (\ref{SDP-2}), {it holds}
\be c^{\top}x^0= -(U^0+W^0_{m_0})\bullet A_0.\label{17-18}\ee

\end{lemma}
{\bf Proof.}
To prove the {lemma}, we will algorithmically construct the number $ m_0$ and the matrices (\ref{G0}).\\

{\it Iteration $\#$ 0.} Consider the following SIP problem:
\be \displaystyle\min_{(x,\mu)}\ \mu, \mbox{ s.t. }
t^{\top}{\cal A}(x)t+\mu\geq 0, \; t \in T,\label{1a}\ee
with the set $T$ defined in (\ref{SetT}). If there exists a feasible solution $(\bar x,\bar \mu)$ of this problem with $\bar \mu<0$, then set
$m_0:=0$   and GO TO the {\it Final step.}

Otherwise for any $x\in X$, the vector $(x,\mu^0=0)$  is an optimal solution of problem (\ref{1a}). It should be noticed
that in problem (\ref{1a}), the index set $T$ is a compact, and the constraints of this problem satisfy the
Slater condition.  Hence, (see e.g. \cite{Bon}), there exist indices and numbers
    \begin{equation*}\tau(i)\in T, \; \; \gamma(i)>0,i \in I_1, \;\; 1\leq |I_1|\leq n+1,\end{equation*}

such that
\be \sum\limits_{i\in I_1}\gamma(i)(\tau(i))^{\top}
A_j\tau(i)=0, \; j=0,1,...,n;\;\;  \sum\limits_{i\in I_1}\gamma(i)=1.\label{2a}\ee
It follows from (\ref{2a}) that the set $I_1$ is nonempty and $\tau(i)\in T_{im}, \,i \in I_1.$

Let us set   $ \ {\beta_1(i):=\sqrt{\gamma(i)}, i\in I_1, \;\;} V^0_1:=(\tau^0(i):={\beta_1(i)} \tau(i), i \in I_1),\; U^0_1:= V^0_1( V^0_1)^{\top}.$

Then equalities (\ref{2a}) take the form
\be U^0_1\bullet A_j=0,\; j=0,1,...,n.\label{A0}\ee
 Denote ${\cal T}_1:={\rm conv} \ \{\tau(i), i \in I_1\}$ and proceed to the next iteration.\\

{\it Iteration $\#$ 1.}
Consider the problem
\begin{equation}\begin{split}& \min\limits_{(x,\mu)}\ \mu, \\
\mbox{ s.t. } \ {\cal A}(x)\tau(i)\geq 0, \; i \in I_1, & \ \ \ t^{\top}{\cal A}(x)t+\mu\geq 0, \; t \in \{t \in T:\rho(t,{\cal T}_1)\geq
\ep_1\},\label{bar3a}\end{split}\end{equation}
where $\ep_1>0$ is  such a number that the set of feasible solutions of problem (\ref{bar3a}) with $\mu=0$ coincides with the set $X$ of feasible
solutions of problem (\ref{SDP-2}). According to Lemma \ref{L0}, such $\ep_1>0$ exists.

If there exists a feasible solution $(\bar x,\bar \mu)$ of problem (\ref{bar3a}) with $\bar \mu<0$,  then STOP and GO TO the {\it Final step} with $m_0:=1.$

Otherwise,  $(x,\mu^0=0)$ with  any  $x\in X$ is an optimal solution of problem (\ref{bar3a}). In the SIP problem (\ref{bar3a}), the index set $\{t \in
T:\rho(t,{\cal T}_1)\geq \ep_1\}$ is compact,
and the constraints  satisfy the following Slater type condition:
$$\exists (\hat x,\hat \mu) \mbox{ such that }  {\cal A}(\hat x)\tau(i)\geq 0, \;
 i \in I_1,\; t^{\top}{\cal A}(\hat x)t+\hat \mu> 0, \; t \in \{t \in T:\rho(t,{\cal T}_1)\geq \ep_1\}.$$

  Hence (see \cite{levin}) there exist indices and numbers
\be \tau(i)\in \{t \in T:\rho(t,{\cal T}_1)\geq \ep_1\}, \;  \gamma(i)>0,i \in \Delta I_1, \;\; 1\leq |\Delta I_1|\leq n+1 ,\label{4.1}\ee

and vectors
$ \lambda^1(i)\in \mathbb R^p_+, \; i\in I_1,$
such that
\be \sum\limits_{i\in \Delta I_1}\gamma(i)(\tau(i))^{\top}A_j\tau(i)+\sum\limits_{i\in I_1}(\lambda^1(i))^{\top}A_j\tau(i)=0, \; j=0,1,...,n; \;\;
\sum\limits_{i\in \Delta I_1}\gamma(i)=1.\label{4a}\ee

It follows from (\ref{4a}) that $\Delta I_1\not=\emptyset$ and for all $x\in X$, it holds
$$(\tau(i))^{\top}{\cal A}(x)\tau(i)=0,\; i \in \Delta I_1;\; (\lambda^1(i))^{\top}{\cal A}(x)\tau(i)=0, \ i \in I_1 .$$
 Hence {$\tau(i)\in T_{im}, \,i \in \Delta I_1.$ From (\ref{2a}) and (\ref{4a}), we get
\be \sum\limits_{i\in  I_2}\gamma(i)(\tau(i))^{\top}A_j\tau(i)+\sum\limits_{i\in I_1}(\lambda^1(i))^{\top}A_j\tau(i)=0, \; j=0,1,...,n,\label{5a}\ee
where $I_2:=I_1\cup\Delta I_1.$ Let us apply to the data set
\be \{\tau(i),\; \gamma(i), \;i\in \Delta I_1;\;\;\; \tau(i),\; \lambda^1(i),\; \gamma(i), i \in I_1\}\label{18-7}\ee
a procedure which is described below.

\begin{itemize}
\item[] \textbf{Procedure} \textbf{\it DAM} (Data Modification).

{\it The Procedure starts with an initial data set}
\be\{\tau(i),\; \gamma(i), \;i\in \Delta I;\;\;\; \tau(i),\; \lambda(i),\; \gamma(i), i \in I\}\label{nach}\ee
{\it such that}
$$ \tau(i)\in T_{im}, \; \tau(i)\not\in {\rm conv}\{\tau(i), i \in I\},\; \gamma(i)>0, \;i\in \Delta I;$$
$$\tau(i)\in T_{im}, \; \lambda(i)\in \mathbb R^p_+ ,\; \gamma(i)>0, i \in I.$$

{\it Set} $P_{+}(i):=\{k \in P:\tau_k(i)>0\},$ $ i \in \Delta I\cup I.$  {\it If}
\be P_{+}(i)\cap (P\setminus P_{+}(s))\not=\emptyset \;\;  \forall  s \in \Delta I \;\;  \forall  i \in I,\label{18-8}\ee
{\it  then STOP. The Procedure \textbf{DAM} is complete.}

{\it If} (\ref{18-8}) {\it  is not satisfied, then find} $s_0\in \Delta I$ {\it and} $i_0\in I$ {\it such that}
\be P_{+}(i_0)\subset P_{+}(s_0).\label{18-9}\ee
{\it Set} $\theta: =\min\limits_{k\in P}\theta_k>0,\;
\mbox{ where }\theta_k:=\left\{\begin{array}{l}
\infty, \ \mbox{if }  k \in P\setminus P_{+}(i_0),\\
\tau_k(s_0)/\tau_k(i_0), \ \mbox{if } k \in P_{+}(i_0).\end{array}\right. $

{\it Let us show that} $\theta <1.$
{\it Suppose the contrary:} $\theta \geq 1.$ {\it Hence} $\theta_k\geq 1 \ \forall k \in P_{+}(i_0),$ {\it and consequently,}
$\tau_k(s_0)\geq \tau_k(i_0)>0,$ $k \in P_{+}(i_0).$
 {\it Notice that since}
$$1=\mathbf e^{\top}\tau(i_0)=
\sum\limits_{k\in P_{+}(i_0)}\tau_k(i_0)\leq \sum\limits_{k\in P_{+}(i_0)}\tau_k(s_0)\leq \sum\limits_{k\in
P_{+}(s_0)}\tau_k(s_0)=\mathbf e^{\top}\tau(s_0)=1,$$
{\it  we conclude that}
$$\sum\limits_{k\in P_{+}(i_0)}\tau_k(s_0)= 1,\; \sum\limits_{k\in P_{+}(i_0)}\tau_k(i_0)= 1,\; \mbox{ and }
\tau_k(s_0)\geq \tau_k(i_0)>0 \ \forall k \in P_{+}(i_0).$$
{\it It follows from the {latter}  conditions   that $\tau(s_0)=\tau(i_0)$
which contradicts the assumption  $\tau(s_0)\not \in conv\{\tau(i), i \in I\}.$
The contradiction proves that $\theta < 1$.

Since, by construction, $\theta$ is strictly positive,
 then the double inequality  $0<\theta<1$ is valid.}

{\it In the  data set (\ref{nach}), let us perform the following replacements:}

\vspace{-6mm}

\begin{equation*}\begin{split} &\tau(s_0)\; \longrightarrow\; \bar \tau(s_0)=(\tau(s_0)-\theta \tau(i_0))/(1-\theta)\geq 0, \; \mathbf e^{\top}\bar \tau(s_0)=1, \;
\bar \tau(s_0)  \not \in {\rm conv}\{\tau(i), i \in I\};\\
&\lambda(i_0)\; \longrightarrow\; \bar \lambda(i_0)=\lambda(i_0)+2\gamma(s_0)\theta(1-\theta)\bar \tau(s_0)\geq 0;\\
&\gamma(i_0)\; \longrightarrow\; \bar \gamma(i_0)=\gamma(i_0)+\gamma(s_0)\theta^2>0;\\
&\gamma(s_0)\; \longrightarrow\; \bar \gamma(s_0)=\gamma(s_0)(1-\theta)^2>0.\end{split}\end{equation*}

\vspace{-2mm}

{\it All other data remain unchanged.}

\vspace{2mm}

{\it For the modified data set, check condition (\ref{18-8}).  If it  is satisfied, then STOP,
the procedure is complete. If (\ref{18-8}) is not satisfied, then  find new indices
$ s_0\in \Delta I$ and $i_0\in I$ such that inclusion (\ref{18-9}) is valid and repeat the steps described above.

The Procedure \textbf{DAM} is completely described.}\\

\end{itemize}

Let us continue proving the Lemma. Recall that we are performing the Iteration $\# 1$ of the algorithm. Having applied the Procedure \textbf{DAM} to the data set (\ref{18-7}), one obtains a new (modified) data set in the  same form  (\ref{18-7})  such that

 \vspace{2mm}

$\bullet$ the indices $\tau(i), i \in I_1,$ are the same as in the initial data set (i.e.,  the procedure  leaved these indices unchanged);

$\bullet$ the modified  indices $\tau(i), $ $ i\in \Delta I_1,$ are the immobile ones in problem (\ref{SDP-1});

$\bullet$ for the modified indices $\tau(i)$ and numbers $\gamma(i),$ $ i \in \Delta I_1,$ relations (\ref{4.1}) are fulfilled;

$\bullet$ for the modified vectors $\lambda^1(i)$ and numbers $\gamma(i),$ $ i \in  I_1,$  it holds
$\lambda^1(i)\in \mathbb R^p_+,$ $ \gamma(i)>0, i \in I_1;$

$\bullet$ for the modified data set (\ref{18-7}), relations (\ref{18-8}) with $\Delta I=\Delta I_1$, $ I=I_1$ and (\ref{5a}) are satisfied.

 \vspace{5mm}

Using the new data (obtained as the result of applying  the Procedure \textbf{DAM}  to the initial data set (\ref{18-7})),
 denote:
$$ \ \beta_2(i):=\sqrt{\gamma(i)},\;i\in I_2,\; V^0_2:=(\beta_2(i)\tau(i), \ i \in I_2),\; L^0_1:=(
\lambda^1(i)/\beta_1(i), \; i \in I_1),$$
$$U^0_2:= V^0_2( V^0_2)^{\top},\; W^0_1:=L^0_1(V^0_1)^{\top},
\; D^0_1:=L^0_1(L^0_1)^{\top}.$$
Then relations (\ref{5a}) can be written as follows:
\be ( U^0_2+  W^0_1)\bullet A_j=0, j=0,1,...,n.\label{A1new}\ee
 GO TO the next iteration.

{\it Iteration  $\#$ $m$, $m\geq 2$.} By the beginning of  the
iteration,  the numbers $\beta_m(i)>0, i \in I_m,$ as well as the  indices, vectors and numbers
$$\tau(i)\in T_{im},\,\,  \gamma(i)>0, \; i \in I_m=I_{m-1}\cup \Delta I_{m-1},\;
\lambda^{m-1}(i)\in \mathbb R^p_+,\; i \in I_{m-1},$$
are found such that

\vspace{2mm}

$\bullet$ relations (\ref{18-8}) with
$\Delta I=\Delta I_{m-1}\not=\emptyset,$ $I=I_{m-1}$ hold;

$\bullet$ the following  equalities are satisfied:
\be \sum\limits_{i\in  I_m}\gamma(i)(\tau(i))^{\top}A_j\tau(i)+\sum\limits_{i\in I_{m-1}}(\lambda^{m-1}(i))^{\top}A_j\tau(i)=0, \; j=0,1,...,n.\label{18-10}\ee

Using these data, matrix
\begin{equation}V^0_m=(\beta_m(i)\tau(i), \ i \in I_m) \;\; \mbox{ with }  \beta_m(i)>0, \; \ i \in I_m,\label{OO1}\end{equation}
 was constructed.

Denote ${\cal T}_m:={\rm conv}\{\tau(i), i \in I_m\}$ and
consider the problem
\begin{equation}\begin{split} & \min\limits_{(x,\mu)\in \mathbb R^{n+1}} \mu, \\
\mbox{ s.t. } \ \ {\cal A}(x)\tau(i)\geq 0, \; i \in I_m, & \
\ t^{\top}{\cal A}(x)+\mu\geq 0, \; t \in \{t \in T:\rho(t,{\cal T}_m)\geq \ep_m\},\label{bar-m}\end{split}\end{equation}
where $\ep_m>0$ is  such a number that the feasible set  of problem (\ref{bar-m}) with $\mu=0$ coincides with the     feasible set $X$ of
  problem (\ref{SDP-2}). According to Lemma \ref{L0}, such $\ep_m$ exists.

If there exists a feasible solution $(\bar x,\bar \mu)$ of problem (\ref{bar-m}) with $\bar \mu<0$, then STOP and GO TO the
  {\it Final step}  with $m_0:=m.$

Otherwise  for any $x\in X,$  vector $(x,\mu^0=0)$  is an optimal solution of problem (\ref{bar-m}).

 Since in problem (\ref{bar-m})
the index set $\{t \in T:\rho(t,{\cal T}_m)\geq \ep_m\}$ is  compact, and the constraints  satisfy the
Slater type condition, then  the optimality of $(x,\mu^0=0)$ provides that there exist indices and numbers
\be \tau(i)\in \mathbb R^p_+, \; \mathbf e^{\top}\tau(i)=1,\;\tau(i)\not \in {\cal T}_m;\;  \gamma(i)>0,i \in \Delta I_m, \;\; 1\leq |\Delta I_m|\leq n
+1,\label{11.1}\ee

and vectors
\be \widehat \lambda^m(i)\in \mathbb R^p_+, \; i\in I_m,\label{4.2*}\ee
that satisfy the following equalities:
\be \sum\limits_{i\in \Delta I_m}\gamma(i)(\tau(i))^{\top}A_j\tau(i)+\sum\limits_{i\in I_m}(\widehat \lambda^m(i))^{\top}A_j\tau(i)=0, \;
j=0,1,...,n;\;\;  \sum\limits_{i\in \Delta I_m}\gamma(i)=1.\label{18-11}\ee

It follows from (\ref{18-11}) that
$$(\tau(i))^{\top}{\cal A}(x)\tau(i)=0,\; i \in \Delta I_m;\; (\widehat \lambda^m(i))^{\top}{\cal A}(x)\tau(i)=0, i \in I_m ,\;\;
\forall\, x\in X,$$
  and, therefore, $\tau(i)\in T_{im},\, i \in \Delta I_m.$

Based on (\ref{18-10}) and (\ref{18-11}), one can conclude that
\be \sum\limits_{i\in  I_{m+1}}\gamma(i)(\tau(i))^{\top}A_j\tau(i)+\sum\limits_{i\in I_m}(\lambda^m(i))^{\top}A_j\tau(i)=0, \; j=0,1,...,n,\label{18-12}\ee
where $I_{m+1}:=I_m\cup\Delta I_m,$ and the vectors $\lambda^m(i)\in \mathbb R^p_+,$ $ i \in I_m,$ are constructed as follows:
$$\lambda^m(i)=\left\{\begin{array}{l}
\lambda^{m-1}(i)+{\widehat \lambda}^m(i), \ i \in I_{m-1},\\
{\widehat \lambda}^m(i), \qquad\qquad \ \ \  i \in \Delta I_{m-1}=I_m\setminus I_{m-1}.\end{array}\right.\qquad\qquad\qquad\qquad\qquad\qquad$$

{Having applied} the described above Procedure \textbf{DAM} to the data set
\be \{\tau(i),\; \gamma(i), \; i\in \Delta I_m;\;\;\; \tau(i),\; \lambda^m(i),\; \gamma(i), i \in I_m \},\label{18-13}\ee
one  will get the modified data set (in the same form) such that

 \vspace{2mm}

$\bullet$ the indices $\tau(i), i \in I_m,$ are the same as in the initial data set
(these indices are not changed  by the Procedure \textbf{DAM});

$\bullet$ the modified indices $\tau(i), $ $ i\in \Delta I_m,$  are the immobile ones in problem (\ref{SDP-1});

$\bullet$ for the modified indices $\tau(i)$ and numbers $\gamma(i),$ $ i \in \Delta I_m,$ relations (\ref{11.1}) are fulfilled;

$\bullet$ for the modified vectors $\lambda^m(i)$ and numbers $\gamma(i),$ $ i \in  I_m,$ it holds
$$\lambda^m(i)\in \mathbb R^p_+,\; \gamma(i)>0,\;\; i \in I_m;$$

$\bullet$ the modified data (\ref{18-13}) satisfy relations (\ref{18-8}) with $\Delta I=\Delta I_m$, $ I=I_m$ and relations (\ref{18-12}).

\vspace{2mm}

Using these new data, let us set
$$\beta_{m+1}(i):=\sqrt{\gamma(i)}, i \in I_{m+1}; \; V^0_{m+1}:=(\beta_{m+1}(i)\tau(i), \ i \in I_{m+1}),\; L^0_m:=(
\lambda^m(i)/\beta_m(i), \; i \in I_m),$$
$$U^0_{m+1}:= V^0_{m+1}( V^0_{m+1})^{\top},\; W^0_m:=L^0_m(V^0_m)^{\top},\;
D^0_m:=L^0_m(L^0_m)^{\top},$$
 where matrix $V^0_m$ was  defined at the previous iteration according to (\ref{OO1}).
Then relations (\ref{18-12}) can be written in the form:
\be ( U^0_{m+1}+  W^0_m)\bullet A_j=0, j=0,1,...,n.\label{A2new}\ee

 Perform the next Iteration $\#(m+1)$.

\vspace{5mm}

{\it Final step.}
It will be proved  in Lemma \ref{final} (see below) that the algorithm consists of a finite number of iterations.

 Hence, for some $0\leq m_0<\infty$,  one of the following  situations will arise:

$\qquad a)$ $m_0=0$ and for problem (\ref{1a}) there exists  a feasible solution $(\bar x, \bar \mu)$ with $\bar \mu<0$;

$\qquad b)$  $m_0>0$ and for  problem (\ref{bar-m})  with $m=m_0$ there exists  a feasible solution $(\bar x, \bar \mu)$ with $\bar \mu<0.$

\vspace{1mm}

 In situation $ a)$   the constraints of the original problem (\ref{SDP-2}) satisfy the Slater condition. Hence, according to the well-known optimality
conditions (see \cite{AH 2013}, for example), if $x^0$ is an optimal
solution of problem (\ref{SDP-2}), then there exists a matrix $U^0\in  \mathcal{CP}^p$ such that
$$U^0\bullet A_j=c_j, \; j=1,2,...,n;\; U^0\bullet \mathcal A(x^0)=0.$$
It follows from the relations above that $U^0$ is a feasible solution of the dual problem (\ref{dual0}) and equality (\ref{17-18}) {holds}.

 Consider situation $ b)$: $m_0>0$.
By the beginning of the final step,  the  matrices
$$ U^0_s,  W^0_s,\; D^0_s,\; s=1,...,m_0-1;\;\;  W_0=\mathbb O_{p},\; U^0_{m_0}=V^0_{m_0}(V^0_{m_0})^\top, $$
 the immobile indices  $\tau(i)$, and numbers $\beta_{m_0}(i)>0,$ $ i \in   I_{m_0},$ have been  constructed.

\vspace{1mm}

Consider the problem
\begin{equation}\begin{split}& \min\limits_{x}c^{\top}x\  , \\
\mbox{ s.t. } {\cal A}(x)\tau(i)\geq 0, \; i \in I_{m_0}, &\ t^{\top}{\cal A}(x)t\geq 0, \; t \in \{t \in T:\rho(t,{\cal T}_{m_0})\geq
\ep_{m_0}\},\label{A1}\end{split}\end{equation}
where  $\ep_{m_0}>0$ is the number used when  problem (\ref{bar-m})  with
$m=m_0$ was formulated.
The way   the number $ \ep_{m_0}$ has been chosen guarantees that the feasible set of problem (\ref{A1})
coincides with the feasible set of problem (\ref{SDP-2}).

Problem (\ref{A1}) satisfies a Slater type condition since, by construction,
$$ {\cal A}(\bar x)\tau(i)\geq 0, \; i \in I_{m_0},\;\; t^{\top}{\cal A}(\bar x)t\geq -\bar \mu>0, \; t \in \{t \in T:\rho(t,{\cal T}_{m_0})\geq \ep_{m_0}\},$$
and the index set $\{t \in T:\rho(t,{\cal T}_{m_0})\geq \ep_{m_0}\}$ is compact.

Let $x^0$ be an optimal solution of problem (\ref{SDP-2}). Then vector $x^0$ is optimal in problem (\ref{A1}) as well. Hence, there exist indices,
numbers and vectors
$$\tau(i)\in \mathbb R^p_+,\; \mathbf e^{\top}\tau(i)=1, \; \tau(i)\not \in {\cal T}_{m_0},\; \gamma(i)>0, i \in I; \ \lambda^{m_0}(i)\in \mathbb R^p_+,\; i \in
I_{m_0},$$
such that
\be \sum\limits_{i\in I}\gamma(i)(\tau(i))^{\top}A_j\tau(i)+\sum\limits_{i\in I_{m_0}}(\lambda^{m_0}(i))^{\top}A_j\tau(i)=c_j, \; j=1,...,n,\label{A3}\ee
\be (\tau(i))^{\top}{\cal A}(x^0)\tau(i)=0,\; i \in I;\; (\lambda^{m_0}(i))^{\top}{\cal A}(x^0)\tau(i)=0,\; i \in I_{m_0}.\label{A5}\ee
Let us set
$$V^0:=(\tau(i)\sqrt{\gamma(i)}, \ i \in I),\;\; L^0_{m_0}:=(\lambda^{m_0}(i)/\beta_{m_0}(i), i \in I_{m_0}),$$
$$U^0:= V^0( V^0)^{\top},\; W^0_{m_0}:=L^0_{m_0}(V^0_{m_0})^{\top},
 D^0_{m_0}:=L^0_{m_0}(L^0_{m_0})^{\top}.$$
 Then relations (\ref{A3}) take the form
\be ( U^0+  W^0_{m_0})\bullet A_j=c_j, j=1,...,n.\label{A6}\ee

It follows from (\ref{A0}), (\ref{A1new}), (\ref{A2new}), and (\ref{A6})
that the constructed set  of matrices (\ref{G0})
is a  feasible  solution of  problem (\ref{dual}).

It was shown above that for any feasible solution $x \in X$  of problem  (\ref{SDP-2}) and any feasible solution (\ref{dGam})
 of problem  (\ref{dual}), the inequality (\ref{17-3}) holds.
   From   the equalities (\ref{17-11-z}) and (\ref{A5}),  it follows that  the feasible solution $x^0\in   X$
of the primal problem (\ref{SDP-2}) and the constructed above feasible solution (\ref{G0})  of problem  (\ref{dual})
turn the inequality (\ref{17-3}) into equality.  The lemma is proved. $\ \ \Box$

\begin{lemma}\label{final} The described in the proof of Lemma \ref{17-L-2} algorithm is finite (i.e. it stops after a finite number of iterations).
\end{lemma}
{\bf Proof.}
If the algorithm has stopped on the Iteration $\#\; 0$ or the Iteration $\#\; 1$, the lemma is proved.
Otherwise let us consider an Iteration $\#\; m$ of the algorithm for some $m\geq 2$.

 At the beginning of this iteration,  we have the set of indices
$\tau(i)\in \mathbb R^p_+, \; i \in I_m,$
where $$I_m=I_{m-1}\cup \Delta I_{m-1}=\Delta I_{0}\cup \Delta I_{1}\cup ... \cup \Delta I_{m-1}
 \;\mbox{ and } \Delta I_{0}:=I_1,\; \Delta I_{s}\not=\emptyset, s=0,...,m-1.$$
As before, denote
$\ P_{+}(i):=\{k \in P: \tau_k(i)>0\}, \ i \in \Delta I_{s}, \; s=0,1,...,m-1.$

Let $i_s$ be an index from the set $\Delta I_{s}$:
$i_s\in \Delta I_{s}, \; \ s =0,1,...,m-1.$

For any $k$, $ 2\leq k\leq m,$ and any $ s$, $ 0\leq s\leq k-2,$
by construction, it holds $$i_s\in I_{k-1}=\Delta I_{0}\cup \Delta I_{1}\cup ... \cup \Delta I_{k-2}, \ i_{k-1}\in \Delta I_{k-1},$$ and
the relations (\ref{18-8}) are fulfilled with $\Delta I=\Delta I_{k-1}$ and  $I=I_{k-1}.$
Hence
$$ P_{+}(i_s)\cap(P\setminus P_{+}(i_{k-1}))\not =\emptyset, \; \  s=0,1,...,k-2,\; k=2,...,m,$$
wherefrom we conclude

\be P_{+}(i_s)\not \subset P_{+}(i_{k-1}), \; \  s=0,1,...,k-2, \; k=2,...,m.\label{25-1}\ee
Consequently, all the sets
$  P_{+}(i_{s}), s=0,1,...,m-1,$ are different.

\vspace{5mm}

Taking into account that on each Iteration $\#$ $m$ it holds $\Delta I_{s}\not=\emptyset,$ $s=0,1,...,m-1,$
one can conclude that the number  $m_0$ of the iterations fulfilled by the algorithm, cannot be greater than some finite number $m_*$,
where $m_*$ is the  maximal  number of all different subsets of the set $P$  satisfying (\ref{25-1}). The lemma is proved. $\qquad \qquad \qquad \qquad \qquad \qquad \qquad \qquad \qquad\qquad \qquad \qquad\Box$

\begin{remark}

The main contribution of the  algorithm used in the proof of Lemma \ref{17-L-2}, consists in the justification of the existence of a finite
number  $m_0$ and the corresponding feasible solution
 (\ref{G0}) of problem  (\ref{dual})
for which equality (\ref{17-18}) is satisfied.

It worth to mention  that it was not the aim of this paper to find a \textquotedblleft good \textquotedblright
 estimate of the minimal value of the number  $m_0.$

Notice also that it is possible that someone can offer  other (pehaps more complex) procedures for finding the finite sets of
matrices  (\ref{G0})  satisfying the constraints of problem (\ref{dual}) and the  equality (\ref{17-18}). Some of such procedures may provide
a better (smaller  than $m_*$) estimate  of the number $m_0.$
\end{remark}

\begin{remark}
In the case of isolated immobile indices, the set of immobile indices is finite:
$T_{im}=\{t^*(j), j \in J_*\}, |J_*|<\infty$ (see Proposition 2.5 in \cite{KT-Opt}). Then on each Iteration $\# \ m$ of the algorithm it holds
$$\Delta I_{m}\not=\emptyset, \ \{\tau(i),i \in \Delta I_{m}\} \subset \{t^*(j), j \in J_*\},$$
$$\{\tau(i),i \in \Delta I_{k}\}\cap \{\tau(i),i \in \Delta I_{s}\}=\emptyset \; \forall k=1,...,m \; \forall s=1,...,m;\; k\not=s,$$ and
 relations (\ref{18-8}) are satisfied. Hence is this case one does not need to use the Procedure  \textbf{DAM} and {has} $m_0\leq |J_*|.$
\end{remark}

The main result of the paper can be formulated in the form of the following theorem which is a consequence of  Lemmas \ref{17-L} and
\ref{17-L-2}.
\begin{theorem}\label{T1} There exists a finite $m_0\geq 0$ such that
problem (\ref{dual}) is dual to the original linear copositive problem {(\ref{SDP-2})} and the strong duality relations are satisfied, i.e.  if the primal problem {(\ref{SDP-2})} admits an optimal solution $x^0$, then the dual problem also has an optimal solution in the form (\ref{G0}) and
equality (\ref{17-18}) holds.
\end{theorem}

\begin{remark}
In our recent paper \cite{KT-SetValued}, we have suggested another strong
 dual formulation  for Linear Copositive Programming. This  formulation
{was} based on the knowledge of the extremal points of the set ${\rm conv}\, T_{im}.$
 In the present paper,  the extended dual problem for  the linear copositive problem (\ref{SDP-2}) is also obtained using the concept  and the properties  of the normalized immobile index set,  but in its final formulation,  we do  not use neither
 the elements of this set (the immobile indices), nor the extremal  points of its convex hull.
\end{remark}

At the end of this section, we would like to  note that as far as we know, with the exception of the mentioned above
  paper \cite{KT-SetValued}, all previously published optimal conditions and duality results for Linear Copositive Programming are formulated
	under  the Slater condition. Here we do not suppose that the Slater condition is satisfied.
    All of this demonstrates the importance and novelty of the results of the paper.

\section{Linear SDP}

Consider a linear SDP problem

\begin{equation} \displaystyle\min_{x}\; c^{\top}x \;\;\;
\mbox{s.t. }  {\cal A}(x)
\in {\cal P}(p). \label{SDP-new}\end{equation}
Following \cite{Ramana-W}, let us adduce the M. Ramana et al.' extended dual for {this} problem:
\begin{eqnarray}&\max\;  -({\widetilde U}+ {\widetilde W}_{m_0})\bullet A_0,\nonumber\\
&\mbox{s.t. \; }( {\widetilde U}_m+ {\widetilde W}_{m-1})\bullet A_j=0,\; j=0,1,...,n,\; m=1,...,m_0,\nonumber
\\&
\mbox{ (\bf{ED-R})}: \qquad\qquad  ({\widetilde U}+{\widetilde W}_{m_0})\bullet A_j=c_j,
 \; j=1,2,...,n; \widetilde U\in {\cal P}(p),\;\; \widetilde W_0=\mathbb O_{p},\qquad\quad
 \nonumber\\
&\qquad \left(\begin{array}{cc}
 \widetilde U_m&  \widetilde W_m\cr
\widetilde W_m^{\top}&  I
\end{array}\right)\in {\cal P}(2p),\; m=1,...,m_0.\label{d1}\end{eqnarray}

It is easy to notice that the  new dual  problem (\ref{dual})
obtained in this paper  for  problem {(\ref{SDP-2})} has a similar structure and properties  as  the dual problem (\textbf{ED-R}) for SDP problem (\ref{SDP-new}).
Nevertheless,  it is worth  mentioning that these  dual problems were obtained using different approaches:
 the dual  problem (\ref{dual}) was formulated and its properties were established using (implicitly)
the concept of the immobile indices while
the dual SDP problem (\textbf{ED-R}) (referred in \cite{Ramana-W} as
the {\it regularized dual problem} (DRP)) was {derived} using the notion of
 the minimal cone which was described there as the output of a special procedure.

To compare these {results}, let us apply the approach, developed in this paper for Linear Copositive Programming, to the SDP problem (\ref{SDP-new}).
Having repeated the described in Section \ref{ExDual} process of building the dual problem, one can  obtain
the extended dual  to problem (\ref{SDP-new})  in the form
\begin{eqnarray} &
\max \ -( U+ W_{m_0})\bullet A_0,\nonumber\\
\mbox{s.t. }&( U_m+ W_{m-1})\bullet  A_j=0,\ j=0,1,...,n,\; m=1,...,m_0,\qquad\quad \nonumber\\
\mbox{(\bf{ED})}:\qquad&( U+ W_{m_0})\bullet  A_j=c_j, \; j=1,2,...,n,\  U\in {\cal P}(p),\;\; W_0=\mathbb O_{p},\nonumber\\
&\left(\begin{array}{cc}
 U_m&  W_m\cr
 W_m^{\top}&  D_m
\end{array}\right)\in {\cal P}(2p),\; m=1,...,m_0.\label{21-1}\end{eqnarray}

The only difference in formulations (\textbf{ED}) and (\textbf{ED-R}) consists {of} the right lower blocks
of the matrices (\ref{d1}) and (\ref{21-1}). Let us show that these problems  are equivalent.

In fact, let $({\widetilde U}_m,\;{\widetilde W}_m, \;  m=1,...,m_0{, \ {\widetilde U}})$ be a feasible solution of problem (\textbf{ED-R}). It is evident that
$ ( U_m={\widetilde U}_m,\;  W_m={\widetilde W}_m,  D_m=I,\;
m=1,...,m_0{, \ U={\widetilde U}})$ is a feasible solution of problem (\textbf{ED}) with the same value of the cost function.

Now let us show that for any feasible solution
\be ( U_m,\;  W_m,  D_m,\; m=1,...,m_0,\; U)\label{fisU}\ee
 of problem (\textbf{ED}) there exists a feasible solution
\be (\;{\widetilde U}_m,\; {\widetilde W}_m, \;
m=1,...,m_0,\; \widetilde U)\label{fistilde}\ee of problem (\textbf{ED-R})  with the same value of the cost function.

Notice that for $m=1,...,m_0, $ it follows from the inclusion (\ref{21-1}) that there exists a matrix
$$ B_m=\left(\begin{array}{c}
 V_m\cr
 {L}_m
\end{array}\right)\in\mathbb R^{2p\times k(m)}\; \mbox{ with}\;\;  V_m\in \mathbb R^{p\times k(m)}\;  \mbox{and}\;  L_m\in \mathbb
R^{p\times k(m)},$$
 such that
$$ {\left(\begin{array}{cc}
 U_m&  W_m\cr
 W_m^{\top}&  D_m
\end{array}\right)}= B_m  B^{\top}_m=\left(\begin{array}{c}
V_m\cr
 {L}_m
\end{array}\right)( V^{\top}_m\; \;  {L}^{\top}_m).$$
Hence, for  $m=1,...,m_0,$ matrices $ U_m,\;  W_m,\; D_m$ admit representations
$$ U_m= V_m V^{\top}_m,\;\;  W_m= V_m L^{\top}_m,\;  D_m= L_m  L^{\top}_m,$$
with some matrices
$ V_m,\;  L_m.$

Let  (\ref{fisU}) {be} a
 feasible solution of problem (\textbf{ED}).
Set
$${\widetilde U}:= U,\; {\widetilde W}_{m_0}:= W_{m_0},\;\; \rho(m_0):=\max\{1, \mu_{max}(
L^{\top}_{m_0}L_{m_0})\},$$
$${\widetilde U}_{m_0}:= \rho(m_0) U_{m_0} ,\; {\widetilde W}_{m_0-1}:= \rho(m_0) W_{m_0-1}.$$
 Here $\mu_{max}(Q)$ denotes  the maximal eigenvalue of matrix $Q\in \mathbb R^{p\times p}.$

It is easy to check that, by construction, we have
$$ ({\widetilde U}+ {\widetilde W}_{m_0})\bullet A_j=c_j, \; j=1,2,...,n;\;\;
({\widetilde U}_{m_0}+ {\widetilde W}_{m_0-1})\bullet A_j=0,\; j=0,1,...,n.$$

Let us show that
\be {\widetilde U}_{m_0}-{\widetilde W}_{m_0} {\widetilde W}^{\top}_{m_0}\in {\cal
P}(p),\label{21-2}\ee or equivalently,
$$t^{\top}\ V_{m_0}(\rho(m_0) I-{\widetilde L}(m_0)) V^{\top}_{m_0}t\geq 0\;\; \forall t \in \mathbb R^p,$$
 or
$${\tau^{\top}(\rho(m_0) I- {\widetilde L}(m_0))\tau\geq 0 \;\; \forall \tau\in \{\tau \in \mathbb R^p: \tau=V^{\top}_{m_0}t,\; t \in
\mathbb R^p\}\subset \mathbb R^p,}$$
where ${\widetilde L}(m_0):= {L}^{\top}_{m_0} {L}_{m_0}.$

It is known (see \cite{Magnus}, p. 230) that for any real symmetric matrix $Q\in \mathcal{S}(p),$ the
inequality $t^{\top}Qt\leq \mu_{max}(Q)t^{\top}t$ is satisfied for any $t\in \mathbb R^p$. Hence
$$\tau^{\top}(\rho(m_0) I-{\widetilde L}(m_0))\tau
=\rho(m_0)\tau^{\top}\tau-
\tau^{\top}{\widetilde L}(m_0)\tau\geq (\rho(m_0)-\mu_{max}({\widetilde L}(m_0)))\tau^{\top}\tau\geq 0
\;\forall \tau\in \mathbb R^p.$$ Inclusion (\ref{21-2}) is proved.

\vspace{5mm}

Suppose that for some $m\leq m_0$ we have constructed matrices
$$ {\widetilde U},\; {\widetilde U}_{m_0},...,{\widetilde U}_{m},\;\;
{\widetilde W}_{m_0},...,{\widetilde W}_{m},\  \widetilde W_{m-1},$$
and such a number $\rho(m)>0$  that $\widetilde W_{m-1}=\rho(m) W_{m-1}$ and  the following relations hold:
$$ (\widetilde U+ \widetilde W_{m_0})\bullet A_j=c_j, \; j=1,2,...,n;$$
$$(\widetilde U_{s}+\widetilde W_{s-1})\bullet A_j=0,\; j=0,1,...,n,\;\; \widetilde U_{s}-\widetilde W_{s}{\widetilde W}^{\top}_{s}\in {\cal P}(p),\;\;s=m_0,m_0-1,...,m. $$

Let us set
$$\rho(m-1):=\max\{1, \rho^2(m)\mu_{max}(L^{\top}_{m-1} L_{m-1})\},$$
$${\widetilde U}_{m-1}:=\rho(m-1) U_{m-1},\;\;  {\widetilde W}_{m-2}:= \rho(m-1) W_{m-2}.$$
 Applying  the described above rules for  the cases where  $m=m_0, m_0-1,...,2,$  we can construct matrices
${\widetilde U},\;{\widetilde U}_{m_0},...,{\widetilde U}_2,$ $ \tilde  W_{m_0},...,{\widetilde W}_2,\; {\widetilde W}_1$
and {the} number $\rho(2)>0$ such that ${\widetilde W}_1=\rho(2) W_1.$

Set $\ \ \rho(1):=\max\{1,\rho^2(2)\mu_{max}({L}^{\top}_{1} {L}_{1})\},\;\; {\widetilde
U}_1:=\rho(1)  U_1.\ $

One can check that the  constructed above matrices  form a feasible solution (\ref{fistilde})  of problem (\textbf{ED-R}) and it holds
$$({\widetilde U}+ {\widetilde W}_{m_0})\bullet A_0=(  U+  W_{m_0})\bullet A_0.$$

Hence, for the SDP problem (\ref{SDP-new}),  we have shown that the dual problem in the form  (\textbf{ED}) is a slight modification of the known dual problem
(\textbf{ED-R}).

\vspace{3mm}

Now, let us compare two  pairs of primal and dual problems:\\

$\qquad\qquad (\alpha)$ the linear copositive problem  {(\ref{SDP-2})} and its dual one (\ref{dual}), and\\

\vspace{-3mm}

$\qquad\qquad (\beta)$ the SDP problem  (\ref{SDP-new}) and  its dual one (\textbf{ED}).

\vspace{3mm}

One can see that these pairs of dual problems are constructed in spaces $S(p)$ and $S(2p)$ using \textbf{the same}  rules,
but their constraints are defined with the help of    \textbf{different}  dual  cones:

\begin{itemize}
\item in the pair  of   problems (\ref{SDP-2}) and (\ref{dual}),
the cone $\mathcal{COP}^p$ is used to formulate the constraints of the
primal copositive problem and the dual cones $\mathcal{CP}^p$  and $\mathcal{CP}^{2p}$
are used to formulate the constraints of the dual one;

\item in the pair  of SDP problems (\ref{SDP-new}) and (\textbf{ED}),
 the cone ${\cal P}(p)$ is used to formulate the constraints of the
 primal SDP problem  and the dual cones ${\cal P}^*(p)={\cal P}(p)$ and ${\cal P}^*(2p)={\cal P}(2p)$ are used for the dual
 formulation.

\end{itemize}

This similarity  points to a deep relationship between these
two classes of conic problems, Linear Copositive Programming and SDP. At the same time, it is  worth mentioning that copositive
problems are more complex and less studied when compared with that of SDP.

\begin{remark} When comparing the complexity of the mentioned above procedures of constructing the pairs of dual problems in SDP and Linear Copositive Programming, notice the following.
\begin{itemize}
\item For SDP problems, one has an  estimate $ m_0\leq \min\{n,p\}$  of the
number $m_0$.  This estimate can be  found using the  fact that the set of immobile indices for an SDP problem is
a subspace of $\mathbb R^p$  and the properties of semi-definite matrices   are well-studied \cite{KT-JMS}.

\item For linear copositive problems, determining  a good estimate of $m_0$ is a much more {challenging} task
as  the set of immobile indices is a union of a  finite number of convex cones
in $ \mathbb R^p$.   Notice that  the  cone  of  copositive matrices and its dual  cone (the cone of completely positive
matrices) are not so well studied (there are many open questions here \cite{Bomze 2012, Dur 2010}).

\item  The cones of copositive and completely positive matrices are   neither self-dual nor homogeneous (see \cite{GOWDA}).
\end{itemize}	
	As it was noticed above,  {finding} a good estimate of the number $m_0$  for copositive problems was not
	our purpose here. 	We plan to  devote a special paper to this  issue.
	\end{remark}

\section{Conclusions and future work}

The main contribution of the paper consists in  developing  a new approach to dual
formulations in Linear Copositive Programming.  This approach permitted us to formulate a new  extended dual problem in explicit form
 and to close the duality gap between the optimal values of the copositive   problem and its extended dual without any CQs or other additional assumptions.

To the best of our knowledge, with the exception of our previous papers \cite{KT-SetValued, KT-Opt2018},
in Linear Copositive Programming, there are  no  other known explicit  strong dual formulations
 that do not require CQs.

In \cite{KT-SetValued, KT-Opt2018}, the  dual  problems were  formulated  based on the explicit  knowledge of
the immobile index set. The  advantage of the dual results presented here if compare with that of the mentioned above results
consists in the fact that now there is no need to find explicitly either the elements of the normalized immobile index set or the extremal points
of its convex hull. For linear copositive problems, the dual formulations obtained in the paper are original and different from that published before.

The new dual formulation for Linear Copositive Programming is similar to the dual formulation for SDP problem proposed by M.Ramana et al.
{\cite{Ramana-W}. This similarity and the fact that the duality results obtained in this paper
{$ \ (i)$} do not use CQs, {$\ (ii)$} have explicit formulation, and {$\ (iii)$} are strong, motivate us to study other applications of the developed approach based the notion of the immobile indices.

In our future work,  we are going to find a better estimate of the number $m_0$ that is essential
for our dual formulation. To obtain this estimate, it will be necessary to
study  new properties of the extended dual problem and its feasible set. We plan also to apply the  results of the paper for
other  classes of  copositive  problems with the aim to develop new   explicit optimality conditions.

\end{document}